\documentclass{amsart}

\usepackage{amsmath,amsthm,amsfonts}

\theoremstyle{plain}
\newtheorem{theorem}{Theorem}[section]

\newtheorem{lemma}{Lemma}[section]

\numberwithin{equation}{section}

\newcommand{\qbin}[2]{\genfrac{[}{]}{0pt}{}{#1}{#2}}
\newcommand{\qbins}[2]{{\textstyle\genfrac{[}{]}{0pt}{}{#1}{#2}}}
\newcommand{\I}{\mathcal{I}}
\newcommand{\Z}{\mathbb{Z}}

\begin{document}

\title{50 Years of Bailey's lemma}
 
\author{S. Ole Warnaar}
 
\address{Instituut voor Theoretische Fysica, Universiteit van Amsterdam,
Valckenierstraat 65, 1018 XE Amsterdam, The Netherlands \\
present address: Department of Mathematics and Statistics,
The University of Melbourne, Vic 3010, Australia}

\maketitle 

\setcounter{section}{-1}

\section{introduction}
Half a century ago, The Proceedings of the London Mathematical Society
published W.~N.~Bailey's influential paper \textit{Identities of the
Rogers--Ramanu\-jan type}~\cite{Bailey49}. The main result therein, which
was inspired by Rogers' second proof of the Rogers--Ramanujan
identities~\cite{Rogers17} (and also \cite{Rogers94,Dyson43,Bailey47}),
is what is now known as Bailey's lemma.
To celebrate the occasion of the lemma's fiftieth birthday
we present a history of Bailey's lemma in 5 chapters
(or rather sections), covering (i) Bailey's work,
(ii) the Bailey chain (iii) the Bailey lattice
(iv) the Bailey lemma in statistical mechanics, and
(v) conjugate Bailey pairs.

Due to size limitations of this paper 
the higher rank
\cite{ML92,LM93,ML95,Milne97,ASW99,Warnaar99} and
trinomial \cite{AB98,Warnaar98,BMP98}
generalizations of the Bailey lemma 
will be treated at the lemma's centennial in 2049.
More extensive reviews of topics (i), (ii) and (iii), can be found in
\cite[Sec. 3]{Andrews85}, \cite{Paule88} and \cite{Bressoud88}, respectively.

\section{The Bailey lemma}\label{secBL}
In an attempt to clarify Rogers' second proof \cite{Rogers17} of the
Rogers--Ramanujan identities, Bailey \cite{Bailey49} was led to the
following simple observation.
\begin{lemma}
If $\alpha=\{\alpha_L\}_{L\geq 0},\dots,\delta=\{\delta_L\}_{L\geq 0}$
are sequences that satisfy
\begin{equation}\label{Bp}
\beta_L=\sum_{r=0}^L \alpha_r u_{L-r}v_{L+r} \qquad \text{and} \qquad
\gamma_L=\sum_{r=L}^{\infty} \delta_r u_{r-L} v_{r+L},
\end{equation}
then
\begin{equation}\label{abcd}
\sum_{L=0}^{\infty} \alpha_L \gamma_L
=\sum_{L=0}^{\infty} \beta_L \delta_L.
\end{equation}
\end{lemma}
The proof is straightforward and merely requires an interchange of sums.
Of course, in the above suitable convergence conditions need to be imposed
to make the definition of $\gamma$ and the interchange of sums meaningful.

The idea behind Bailey's lemma is clear. When trying to prove a
complicated identity of the form $\sum_L A_L=\sum_L B_L$ it is a considerable
step in the right direction if one can find a dissection of this identity into
two identities of the type \eqref{Bp} where $A_L=\alpha_L \gamma_L$ and
$B_L=\beta_L \delta_L$. Or, as Slater put it in Bailey's 
obituary~\cite{Slater62},
\begin{quote}
\textit{The root of the underlying idea \dots is that of transforming
a doubly infinite series into a simply infinite and a finite series.
In a geometric sense, this involves summing over a triangle instead of over a
rectangle.}
\end{quote}

In applications of his transform, Bailey chose $u_L=1/(q)_L$ and
$v_L=1/(aq)_L$, with the usual definition of the $q$-shifted factorial,
$(a)_{\infty}=(a;q)_{\infty}=\prod_{k=0}^{\infty}(1-aq^k)$ and
$(a)_L=(a;q)_L=(a)_{\infty}/(aq^L)_{\infty}$
for $L\in\Z$. (Throughout, we assume that $0<|q|<1$.)
With this choice, equation \eqref{Bp} becomes
\begin{equation}\label{BP}
\beta_L=\sum_{r=0}^L \frac{\alpha_r}{(q)_{L-r}(aq)_{L+r}} 
\qquad \text{and} \qquad
\gamma_L=\sum_{r=L}^{\infty} \frac{\delta_r}{(q)_{r-L}(aq)_{r+L}}.
\end{equation}
A pair of sequences that satisfies the first equation of \eqref{BP} 
is called a Bailey pair relative to $a$. Similarly, 
the second equation defines a conjugate Bailey pair relative to $a$.

Still following Bailey, one can employ the $q$-Saalsch\"utz summation
\cite[Eq. (II.12)]{GR90} to establish that $(\gamma,\delta)$ with
\begin{equation}\label{gd}
\begin{split}
\gamma_L&=\frac{(\rho_1)_L(\rho_2)_L(aq/\rho_1\rho_2)^L}
{(aq/\rho_1)_L(aq/\rho_2)_L}
\frac{1}{(q)_{M-L}(aq)_{M+L}} \\[2mm]
\delta_L&=\frac{(\rho_1)_L(\rho_2)_L(aq/\rho_1\rho_2)^L}
{(aq/\rho_1)_M(aq/\rho_2)_M}
\frac{(aq/\rho_1\rho_2)_{M-L}}{(q)_{M-L}}
\end{split}
\end{equation}
provides a conjugate Bailey pair.

As we shall see in the next section this conjugate Bailey pair leads
to the very important concept of the Bailey chain.
However, Bailey missed an opportunity here and made the 
(mis)judgement~\cite[Page 4]{Bailey49}:
\begin{quote}
\textit{
These values of $\delta_L,\gamma_L$ \dots lead to \dots 
results involving only terminating basic series.
We are, however, more concerned with identities of the Rogers--Ramanujan type
in this paper, as the most general formulae for basic series are too
involved to be of any great interest.}
\end{quote}
Consequently, Bailey only considered the conjugate Bailey pair \eqref{gd}
when the parameter $M$ tends to infinity.
Also taking the limit $\rho_1,\rho_2\to\infty$ yields
\begin{equation}\label{gdinf}
\gamma_L=\frac{a^L q^{L^2}}{(aq)_{\infty}} \qquad \text{and} \qquad
\delta_L=a^L q^{L^2},
\end{equation}
which substituted into \eqref{abcd} gives
\begin{equation}\label{ab}
\frac{1}{(aq)_{\infty}}\sum_{L=0}^{\infty} a^L q^{L^2} \alpha_L=
\sum_{L=0}^{\infty} a^L q^{L^2} \beta_L.
\end{equation}
The proof of the Rogers--Ramanujan and similar such identities
requires the input of suitable Bailey pairs into \eqref{ab}.
For example, from Rogers' work \cite{Rogers17} one can infer the
following Bailey pair relative to $1$: $\alpha_0=1$ and
\begin{equation}\label{BPRogers}
\alpha_L=(-1)^L q^{L(3L-1)/2}(1+q^L),\qquad
\beta_L=\frac{1}{(q)_L}.
\end{equation}
Thus one finds
\begin{equation*}
\frac{1}{(q)_{\infty}}\sum_{L=-\infty}^{\infty}
(-1)^L q^{L(5L-1)/2}=\sum_{n=0}^{\infty} \frac{q^{n^2}}{(q)_n}.
\end{equation*}
The application of the Jacobi triple product identity \cite[Eq. (II.28)]{GR90}
yields the first Rogers--Ramanujan identity~\cite{Rogers94,Rogers17}
\begin{equation}\label{RR1}
\sum_{n=0}^{\infty} \frac{q^{n^2}}{(q)_n}=
\frac{1}{(q,q^4;q^5)_{\infty}},
\end{equation}
with the notation $(a_1,\dots,a_k;q)_n=(a_1;q)_n\dots (a_k;q)_n.$
The second Rogers--Ram\-an\-ujan identity
\begin{equation}\label{RR2}
\sum_{n=0}^{\infty} \frac{q^{n(n+1)}}{(q)_n}=
\frac{1}{(q^2,q^3;q^5)_{\infty}}
\end{equation}
follows in a similar fashion using the Bailey pair \cite{Rogers17}
\begin{equation}\label{BPRogers2}
\alpha_L=(-1)^L q^{L(3L+1)/2}(1-q^{2L+1})/(1-q),\qquad
\beta_L=\frac{1}{(q)_L}
\end{equation}
relative to $q$.
By collecting a list of 96 Bailey pairs, 
and using \eqref{ab} or the identity obtained from \eqref{gd}
and \eqref{abcd} by taking $M,\rho_1\to\infty$ and $\rho_2=-q^{k/2}$
(with $k$ a small nonnegative integer)
Slater compiled her famous list of 130  identities of the
Rogers--Ramanujan type~\cite{Slater51,Slater52}. 
The next two sections deal with more systematic ways of finding Bailey pairs.

\section{The Bailey chain}

By dismissing the conjugate Bailey pair \eqref{gd} in its finite form
(i.e., with $M$ finite, or, equivalently, with $\rho_1$ or $\rho_2$
of the form $q^{-N}$) Bailey missed a very
effective mechanism for generating Bailey pairs. Namely, if we substitute
the conjugate pair \eqref{gd} into \eqref{abcd} the resulting equation has
the same form as the defining relation \eqref{BP} of a
Bailey pair. This is formalized in the following theorem due to
Andrews~\cite{Andrews84,Andrews85}.
\begin{theorem}\label{AB}
Let $(\alpha,\beta)$ form a Bailey pair relative to $a$. Then
so does $(\alpha',\beta')$ with
\begin{equation}\label{iter}
\begin{split}
\alpha_L'&=\frac{(\rho_1)_L(\rho_2)_L(aq/\rho_1\rho_2)^L}
{(aq/\rho_1)_L(aq/\rho_2)_L} \alpha_L \\[2mm]
\beta_L'&=\sum_{r=0}^L \frac{(\rho_1)_r(\rho_2)_r(aq/\rho_1\rho_2)^r
(aq/\rho_1\rho_2)_{L-r}}{(aq/\rho_1)_L(aq/\rho_2)_L(q)_{L-r}}\beta_r.
\end{split}
\end{equation}
\end{theorem}
Again letting $\rho_1,\rho_2$ tend to infinity leads to the important special
case
\begin{equation}\label{iter2}
\alpha_L'=a^L q^{L^2}\alpha_L \qquad \text{and}\qquad
\beta_L'=\sum_{r=0}^L \frac{a^r q^{r^2}}{(q)_{L-r}}\beta_r,
\end{equation}
which, for $a=1$ and $a=q$, was also discovered by Paule~\cite{Paule85}.
One now finds that the Bailey pairs \eqref{BPRogers} and \eqref{BPRogers2}
of Rogers can be obtained from the $a=1$ and $a=q$ instances of the Bailey
pair~\cite{Andrews84}
\begin{equation}\label{initial}
\alpha_L=(-1)^L q^{\binom{L}{2}}\frac{(1-aq^{2L})(a)_L}{(1-a)(q)_{L}},
\qquad \beta_L=\delta_{L,0}
\end{equation}
by application of \eqref{iter2}.
The Bailey pair \eqref{initial} is an immediate consequence of the
inverse Bailey transform~\cite{Andrews84}
\begin{equation}\label{inverseBP}
\alpha_L=\frac{1-aq^{2L}}{1-a}\sum_{r=0}^L
\frac{(-1)^{L-r}q^{\binom{L-r}{2}}(a)_{L+r}}{(q)_{L-r}}\beta_r,
\end{equation}
which follows from \eqref{BP} and \cite[Eq. (II.21)]{GR90} specialized to
$aq=bc$.
The iteration of \eqref{iter} or \eqref{iter2} leads to what is known
as the Bailey chain~\cite{Andrews84,Paule88}:
\begin{equation*}
(\alpha,\beta) \to (\alpha',\beta')
\to (\alpha'',\beta'') \to \cdots
\end{equation*}
and thus, given a single Bailey pair, one immediately finds an infinite
sequence of Bailey pairs. (To be compared with the 96 Bailey pairs
collected by Slater!)
As an example, iteration of \eqref{inverseBP}
gives the Bailey pair 
\begin{align*}
\alpha_L&=
(-1)^L a^{kL}q^{kL^2+\binom{L}{2}}\frac{(1-aq^{2L})(a)_L}{(1-a)(q)_L}\\
\beta_L&=\sum_{L\geq n_1\geq\dots\geq n_{k-1}\geq 0}
\frac{a^{n_1+\cdots+n_{k-1}}q^{n_1^2+\cdots+n_{k-1}^2}}
{(q)_{L-n_1}(q)_{n_1-n_2}\cdots(q)_{n_{k-2}-n_{k-1}}(q)_{n_{k-1}}}.
\end{align*}
Substituting this into the defining relation \eqref{BP} of a Bailey pair
and letting $L$ tend to infinity gives, for $a=1$ or $a=q$,
\begin{equation*}
\sum_{n_1,\dots,n_{k-1}}
\frac{a^{n_1+\cdots+n_{k-1}}q^{n_1^2+\cdots+n_{k-1}^2}}
{(q)_{n_1-n_2}\cdots(q)_{n_{k-2}-n_{k-1}}(q)_{n_{k-1}}}
=\frac{1}{(q)_{\infty}}
\sum_{r=-\infty}^{\infty}(-1)^r a^{kr}q^{kr^2+\binom{r}{2}}.
\end{equation*}
Using Jacobi's triple product identity finally yields
\begin{equation}\label{AG}
\sum_{n_1,\dots,n_{k-1}}
\frac{q^{n_1^2+\cdots+n_{k-1}^2+n_i+\cdots+n_{k-1}}}
{(q)_{n_1-n_2} \ldots (q)_{n_{k-2}-n_{k-1}}(q)_{n_{k-1}}}
=\frac{(q^i,q^{2k+1-i},q^{2k+1};q^{2k+1})_{\infty}}
{(q)_{\infty}},
\end{equation}
where $i=1$ or $i=k$.
For $k=2$ these are the Rogers--Ramanujan identities \eqref{RR1} and
\eqref{RR2}, whereas for $k\geq 3$ they are Andrews' analytic counterpart
of Gordon's partition theorem~\cite{Andrews74}.
In fact, Andrews' identities are \eqref{AG} for all $i=1,\dots,k$
and one concludes that the Bailey chain mechanism has failed to
produce all of these. What is required is an extension of the
Bailey chain known as the Bailey lattice. This will be our next topic.
(Prior to the invention of the Bailey lattice
Paule~\cite{Paule85} already obtained \eqref{AG} for all $i$ using
``ad hoc'' Bailey lattice-like transformations.)

\section{The Bailey lattice}
One of the features of Theorem~\ref{AB} is that it transforms a
Bailey pair relative to $a$ into a new Bailey pair relative to $a$.
More generally one can of course try to transform a Bailey pair
relative to $a$ into a Bailey pair relative to $b$. Agarwal, Andrews and
Bressoud have formulated this problem in a general setting of infinite
dimensional matrices~\cite{AAB87,Bressoud88}. Here we shall only be concerned 
with concrete examples of such ``Bailey lattice'' transformations.
Since the parameter $a$ is no longer fixed we shall write 
$(\alpha(a),\beta(a))$ for a Bailey pair relative to $a$.
\begin{theorem}\label{lattice1}
Fix $N$ a nonnegative integer and set $b=aq^N$.
Let $(\alpha(b),\beta(b))$ be a Bailey pair. Then
so is $(\alpha'(a),\beta'(a))$ with
\begin{align*}
\alpha'_L(a) &=(1-aq^{2L})(aq)_N
\frac{(\rho_1)_L(\rho_2)_L(aq/\rho_1\rho_2)^L}
{(aq/\rho_1)_L(aq/\rho_2)_L} \\
&\qquad \times\sum_{j=0}^N (-1)^j a^j q^{2Lj-j(j+1)/2}\qbin{N}{j}
\frac{(aq)_{2L-j-1}}{(aq)_{2L-j+N}}\:\alpha_{L-j}(b) \\
\beta'_L(a) &=\sum_{r=0}^L
\frac{(\rho_1)_r(\rho_2)_r (aq/\rho_1\rho_2)^r(aq/\rho_1 \rho_2)_{L-r}}
{(aq/\rho_1)_L(aq/\rho_2)_L (q)_{L-r}}\:\beta_r(b).
\end{align*}
\end{theorem}
Here we have used the $q$-binomial coefficient defined as
$\qbins{a}{b}=\frac{(q^{a-b+1})_b}{(q)_b}$
for $b\geq 0$ and $0$ otherwise.
A very similar result can be stated as follows.
\begin{theorem}\label{lattice2}
Fix $N$ a nonnegative integer and set $b=aq^N$.
Let $(\alpha(b),\beta(b))$ be a Bailey pair. Then
so is $(\alpha'(a),\beta'(a))$ with
\begin{align*}
\alpha'_L(a) &=(1-aq^{2L})(aq)_N
\sum_{j=0}^N\frac{(\rho_1)_{L-j}(\rho_2)_{L-j}(bq/\rho_1\rho_2)^{L-j}}
{(bq/\rho_1)_{L-j}(bq/\rho_2)_{L-j}} \\
&\qquad \times (-1)^j a^j q^{2Lj-j(j+1)/2}\qbin{N}{j}
\frac{(aq)_{2L-j-1}}{(aq)_{2L-j+N}}\:\alpha_{L-j}(b) \\
\beta'_L(a) &=\sum_{r=0}^L
\frac{(\rho_1)_r(\rho_2)_r(bq/\rho_1\rho_2)^r
(bq/\rho_1\rho_2)_{L-r}}
{(bq/\rho_1)_L(bq/\rho_2)_L(q)_{L-r}}\:\beta_r(b).
\end{align*}
\end{theorem}
The $N=0$ and $N=1$ cases of the first theorem correspond to the
Bailey chain of Theorem~\ref{AB} and the Bailey lattice of
\cite[Lemma 1.2]{AAB87}, respectively.
The second theorem for $N=0$ is again the Bailey chain whereas
for $N=1$ it is a variation of the Bailey lattice of \cite[Lemma 4.3]{SW98}.
Theorem~\ref{lattice1} was also found by Krattenthaler and Foda~\cite{KF}.

First we prove  Theorem~\ref{lattice1}. 
Substituting the expression for $\alpha'(a)$
in the ``primed'' version of \eqref{BP},
transforming $j\to r-j$ and then interchanging
the order of summation, gives
\begin{align*}
\beta'_L(a)&=\sum_{j=0}^L
\frac{(\rho_1)_j(\rho_2)_j(aq/\rho_1\rho_2)^j\alpha_j(b)}
{(aq/\rho_1)_j(aq/\rho_2)_j(q)_{L-j} (bq)_{2j}(aq^{2j+1})_{L-j}} \\
& \times \lim_{a_4\to\infty}{_8W_7}(aq^{2j};a_4,a/b,\rho_1 q^j,\rho_2 q^j,
q^{-L+j};q,abq^{L+j+2}/a_4\rho_1\rho_2),
\end{align*}
where we employed the conventional short-hand notation for very-well-poised
basic hypergeometric series~\cite{GR90}.
By Watson's ${_8\phi_7}$ transformation \cite[Eq. (III.18)]{GR90}
(with $a\to aq^{2j},b\to a_4,c\to a/b,d\to\rho_1q^j,e\to\rho_2q^j$ and 
$n\to L-j$) this can be simplified to
\begin{align*}
\beta'_L(a)&=\sum_{j=0}^L\sum_{r=0}^{L-j}
\frac{(\rho_1)_{j+r}(\rho_2)_{j+r}(aq/\rho_1\rho_2)^{j+r}
(aq/\rho_1\rho_2)_{L-j-r}\alpha_j(b)}
{(aq/\rho_1)_L(aq/\rho_2)_L(q)_{L-j-r}(q)_r(bq)_{r+2j}}.
\end{align*}
Shifting $r\to r-j$, then interchanging sums
and recalling the definition of $\beta_r(b)$, this indeed yields the
second transformation claimed in the theorem.

The second proof proceeds in a similar manner.
Substituting the expression for $\alpha'(a)$
in the ``primed'' version of \eqref{BP},
transforming $j\to r-j$ and then interchanging
the order of summation, gives
\begin{align*}
\beta'_L(a)&=\sum_{j=0}^L 
\frac{(\rho_1)_j(\rho_2)_j(bq/\rho_1\rho_2)^j\alpha_j(b)}
{(bq/\rho_1)_j(bq/\rho_2)_j(q)_{L-j}(bq)_{2j}(aq^{2j+1})_{L-j}}\\
& \qquad \qquad \times \lim_{a_4\to\infty}
{_6W_5}(aq^{2j},a_4,a/b,q^{-L+j};q,bq^{L+j+2}/a_4).
\end{align*}
By Rogers' ${_6\phi_5}$ summation \cite[Eq. (II.21)]{GR90}
(with $a\to aq^{2j},b\to a_4,c\to a/b$ and
$n\to L-j$) this can be simplified to
\begin{equation*}
\beta'_L(a)=\sum_{j=0}^L
\frac{(\rho_1)_j(\rho_2)_j(bq/\rho_1\rho_2)^j\alpha_j(b)}
{(bq/\rho_1)_j(bq/\rho_2)_j(q)_{L-j}(bq)_{L+j}}.
\end{equation*}
Using the $q$-Saalsch\"utz sum \cite[Eq. (II.12)]{GR90}
(with $a\to\rho_1q^j,b\to\rho_2q^j,c\to bq^{2j+1}$ and $n\to L-j$ ) 
this can be rewritten as
\begin{equation*}
\beta'_L(a)=\sum_{j=0}^L\sum_{r=0}^{L-j}
\frac{(\rho_1)_{j+r}(\rho_2)_{j+r}(bq/\rho_1\rho_2)^{j+r}
(bq/\rho_1\rho_2)_{L-j-r}\alpha_j(b)}
{(bq/\rho_1)_L(bq/\rho_2)_L(q)_{L-j-r}(q)_r(bq)_{r+2j}}.
\end{equation*}
Shifting $r\to r-j$, interchanging sums and recalling the definition
\eqref{BP} yields the second expression of Theorem~\ref{lattice2}.

To see that we are now in the position to prove \eqref{AG} for
all $i=1,\dots,k$ we follow \cite{AAB87} and take the Bailey pair
of equation \eqref{initial} with $a=q$ as starting point.
Applying \eqref{iter2} $k-i+1$ times, then Theorem~\ref{lattice1} with $N=1$
and $\rho_1,\rho_2\to\infty$ once, and then again \eqref{iter2} $i-2$ times,
one finds the Bailey pair $\alpha_0=1$,
\begin{align*}
\alpha_L&=(-1)^L q^{kL^2+\binom{L}{2}+(k-i+1)L}(1+q^{(2i-2k-1)L})  \\
\beta_L&=\sum_{L\geq n_1\geq\dots\geq n_{k-1}\geq 0}
\frac{q^{n_1^2+\cdots+n_{k-1}^2+n_i+\cdots+n_{k-1}}}
{(q)_{L-n_1}(q)_{n_1-n_2}\cdots(q)_{n_{k-2}-n_{k-1}}(q)_{n_{k-1}}},
\end{align*}
relative to $1$. Substituting this result into \eqref{BP}, letting $L$ tend
to infinity and using the triple product identity one arrives at the
identities \eqref{AG} for $i=2,\dots,k$.

\section{The Bailey lemma in statistical mechanics}
In section~\ref{secBL} we already mentioned Slater's famous list of
130 identities of Rogers--Ramanujan type~\cite{Slater51,Slater52}. 
She found these identities by exploiting extensive lists of Bailey
pairs (grouped from A to M) extracted from Rogers' or Bailey's papers
\cite{Rogers17,Bailey49} or from known basic hypergeometric function 
identities. For example, the first group of Bailey pairs (all due to Rogers)
reads $\alpha_0=1$,
\begin{equation*}
\begin{tabular}{c|c|c|c}
& $\beta_L$ & $\alpha_{3L\pm 1}$ & $\alpha_{3L}$ \\ \hline &&&\\[-3.5mm]
A(1) & $1/(q)_{2L}$ & $-q^{(2L\pm 1)(3L\pm 1)}$ & $q^{L(6L-1)}+q^{L(6L+1)}$ \\
A(3) & $q^L/(q)_{2L}$ & $-q^{2L(3L\pm 1)}$ & $q^{2L(3L-1)}+q^{2L(3L+1)}$ \\
A(5) & $q^{L^2}/(q)_{2L}$ & $-q^{L(3L\pm 1)}$ & $q^{L(3L-1)}+q^{L(3L+1)}$ \\
A(7) & $q^{L^2-L}/(q)_{2L}$ & $-q^{(L\pm 1)(3L\pm 1)}$ & 
$q^{L(3L-2)}+q^{L(3L+2)}$
\end{tabular}
\end{equation*}
with $a=1$, and
\begin{equation*}
\begin{tabular}{c|c|c|c}
& $\beta_L$ & $\alpha_{3L-(1\mp 1)/2}$ & $\alpha_{3L+1}$  \\
\hline &&&\\[-3.5mm] 
A(2) & $1/(q^2)_{2L}$ & $q^{L(6L\pm 1)}$ & 
$-q^{(2L+1)(3L+1)}-q^{(2L+1)(3L+2)}$ \\
A(4) & $q^L/(q^2)_{2L}$ & $q^{2L(3L\pm 2)}$ & 
$-q^{2(L+1)(3L+1)}-q^{2L(3L+2)}$ \\
A(6) & $q^{L^2}/(q^2)_{2L}$ & $q^{L(3L\mp 1)}$ & $-q^{L(3L+1)}-q^{(L+1)(3L+2)}$ \\
A(8) & $q^{L^2+L}/(q^2)_{2L}$ & $q^{L(3L\pm 2)}$ & 
$-q^{(L+1)(3L+1)}-q^{L(3L+2)}$ 
\end{tabular}
\end{equation*}
with $a=q$. The Bailey pairs given in equations \eqref{BPRogers} and
\eqref{BPRogers2}, which were used by Rogers
to prove the Rogers--Ramanujan identities, are items labelled B(1) and B(3).

Remarkably, in recent work on exactly solvable lattice models in statistical
mechanics identities have arisen 
(see \cite{Berkovich94,FQ95,BM96,Warnaar96,BMS98} and references therein),
which for each pair of integers $(p,p')$,
with $1<p<p'$ and $\gcd(p,p')=1$ imply a family of Bailey pairs~\cite{FQ96}.
Moreover, many of the Bailey pairs of 
Rogers and Slater (as well as later pairs found in 
\cite{Andrews84,Paule85,AAB87}) are included as special cases.

First we need a class of polynomials known as the one-dimensional
configuration sums of the Andrews--Baxter--Forrester model~\cite{ABF84,FB85}.
(See \cite{ABBBFV87} for a partition theoretic interpretation of the
configuration sums.)
For coprime integers $p,p'$ with $1\leq p<p'$, and integers
$1\leq b,s\leq p'-1$, $0\leq r\leq p-1$ and $L\geq 0$ such that 
$L+s+b$ is even, define
\begin{equation*}
X^{(p,p')}_{r,s}(L,b)
=\sum_{j\in\Z} \Bigr\{q^{j(p p'j+p'r-ps)}
\qbins{L}{\frac{L+s-b}{2}-p'j}-
q^{(pj+r)(p'j+s)}\qbins{L}{\frac{L-s-b}{2}-p'j}\Bigl\}.
\end{equation*}
For $r=b-\lfloor (b+1)(p'-p)/p'\rfloor$, (with $\lfloor x \rfloor$ the
integer part of $x$)
the one-dimensional configuration sums are generating functions of
certain sets of restricted lattice paths, and hence are polynomials
with positive coefficients. This is not at all manifest from the
above definition, and the identities referred to in the above claim a
different, manifestly positive representation for the configuration sums.
The simplest of these identities arise when $|p'r-ps|=1$ and $b=s$
which we assume throughout the remainder of this section.

For the moment also assume that $p<p'<2p$, and define nonnegative integers
$\nu_0,\dots,\nu_n$ by the continued fraction expansion
$p/(p'-p)=[\nu_0,\nu_1,\dots,\nu_n]$.
The integers $n$ and $\nu_j$, can be used to further define
$t_m=\sum_{j=0}^{m-1}\nu_j~(1\leq m\leq n)$ and 
$d=-2+\sum_{j=0}^n \nu_j$.
These latter numbers define a so-called
fractional incidence matrix $\I$ and fractional
Cartan-type matrix $B=2I-\I$
(with $I$ the $d$ by $d$ unit matrix) as follows
\begin{equation*}
\I_{i,j} = \begin{cases}
\delta_{i,j+1} + \delta_{i,j-1} & \text{for
$1 \leq i<d,~ i \neq t_m$}, \\
\delta_{i,j+1} + \delta_{i,j} - \delta_{i,j-1} & \text{for
$i=t_m,~1\leq m\leq n-\delta_{\nu_n,2}$}, \\
\delta_{i,j+1} + \delta_{\nu_n,2} \delta_{i,j} & \text{for $i=d$}.
\end{cases}
\end{equation*}
When $p'=p+1$ the matrix $\I$ has entries
$\I_{i,j}=\delta_{|i-j|,1}$ ($i,j=1,\dots,p-2$), so that $B$ corresponds
to the Cartan matrix of the Lie algebra A$_{p-2}$.
When $p=2k-1$ and $p'=2k+1$ one finds
$\I_{i,j}=\delta_{|i-j|,1}+\delta_{i,j}\delta_{i,k-1}$ ($i,j=1,\dots,k-1$),
so that $B$ corresponds to the Cartan-type matrix of the tadpole graph of
$k-1$ nodes.

Using the above definitions we have the following 
result~\cite{Berkovich94,FQ95,BM96,Warnaar96,BMS98}:
\begin{theorem}\label{thmF}
Let $1<p<p'<2p$ with $\gcd(p,p')=1$ and let $r(\leq p'-1)$ and $s(\leq p'-1)$
satisfy $|p'r-ps|=1$. Then
\begin{equation}\label{F}
X_{r,s}^{(p,p')}(2L,s)=\sum_{m\in 2\Z^d}
q^{\frac{1}{4}mBm}\prod_{j=1}^d\qbin{L\delta_{j,1}+\frac{1}{2}(\I m)_j}{m_j}.
\end{equation}
\end{theorem}
Here we use the notation $mBm=\sum_{j,k} m_j B_{j,k}m_k$ and
$(\I m)_j=\sum_k \I_{j,k}m_k$.
The corresponding identities for $2p<p'$ follow simply from the symmetry
$$X_{r,s}^{(p,p')}(L,b;q)=
q^{\frac{1}{4}(L^2-(b-s)^2)} X_{b-r,s}^{(p'-p,p')}(L,b;1/q).$$

Foda and Quano~\cite{FQ96} used (special cases of)
the above theorem and symmetry relation together with the Bailey 
lemma to prove conjectured $q$-series identities for Virasoro characters.
Indeed, we can readily extract the following Bailey pairs
relative to 1~\cite{FQ96,BMSW97}: $\alpha_0=1$,
\begin{equation}\label{alphabeta}
\begin{split}
\alpha_L&=\begin{cases}
q^{j(jpp'+rp'-sp)}+q^{j(jpp'-rp'+sp)}&\text{ for $L=jp'>0$}\\
-q^{(jp\pm r)(jp'\pm s)}&\text{ for $L=jp'\pm s>0$}\\
0 & \text{ otherwise}
\end{cases}\\
\beta_L&=X^{(p,p')}_{r,s}(2L,s)/(q)_{2L},
\end{split}
\end{equation}
where in the expression for $\beta$ the representation \eqref{F}
of $X^{(p,p')}_{r,s}$ is taken.
We note that $(p,p')=(2,3)$ (so that $r=s=1$) corresponds to the Bailey 
pair A(1) and $(p,p')=(1,3)$ ($r=0$ and $s=1$) to A(5).
We also remark that $(p,p')=(2,5)$ ($r=1$, $s=2$) is the Bailey pair
\cite[Eq. (5.3)]{Andrews84}.

\section{Conjugate Bailey pairs}
We have just seen that each pair of coprime integers $(p,p')$ 
labels a Bailey pair. Next we discuss some recent developments
which show that a similar result holds for conjugate Bailey 
pairs~\cite{SW97,SW98,SW99a}.

First we need to introduce the string functions associated to admissible
representations of the affine algebra A$_1^{(1)}$~\cite{KW88}.
Again fix a pair of positive, coprime integers $(p,p')$.
Let $0\leq\ell\leq p'-2$ and let $\Lambda_0,\Lambda_1$ denote the fundamental
weights 
of A$_1^{(1)}$. Then Kac and Wakimoto showed that the A$_1^{(1)}$ character
of the admissible highest weight module of highest weight 
$(p'/p-\ell-2)\Lambda_0+\ell\Lambda_1$ is given by a generalized Weyl-Kac
formula as follows
\begin{equation*}
\chi_{\ell}^{(p,p')}(z,q)=\frac{\sum_{\sigma=\pm 1}\sigma
\Theta_{\sigma(\ell+1),p'}(z,q^p)}
{\sum_{\sigma=\pm 1}\sigma\Theta_{\sigma,2}(z,q)}.
\end{equation*}
Here $\Theta_{n,m}$ is a classical theta function,
$\Theta_{n,m}(z,q)=\sum_{j\in\Z+n/2m} q^{mj^2}z^{-mj}$.
Note that when $p>1$ we are dealing with nonintegral highest weights.
The (normalized) A$_1^{(1)}$ string functions of 
level $p'/p-2$ are defined by the
expansion
\begin{equation*}
\chi_{\ell}^{(p,p')}(z,q)=q^{\frac{1}{8}-\frac{(\ell+1)^2p}{4p'}}
\sum_{m\in\Z}C_{m,\ell}^{(p,p')}(q)z^{-\frac{1}{2}m},
\end{equation*}
which immediately implies that $C_{m,\ell}^{(p,p')}(q)=0$ unless $m+\ell$ is
even.
An explicit expression for the string functions can be derived 
as a double sum of Hecke indefinite modular form type~\cite{KP84,ACT91,SW99a}
\begin{align*}
C_{m,\ell}^{(p,p')}(q)&=\frac{1}{(q)_{\infty}^3}
\Bigl\{\sum_{\substack{i\geq 0\\[.5mm] j\geq 0}}-
\sum_{\substack{i<0 \\[.5mm]j<0}}\Bigr\}
(-1)^i q^{\frac{1}{2}i(i+m)+p'j(pj+i)+\frac{1}{2}(\ell+1)(2pj+i)} \\
&-\frac{1}{(q)_{\infty}^3}
\Bigl\{\sum_{\substack{i\geq 0 \\[.5mm]j>0}}-
\sum_{\substack{i<0\\[.5mm] j\leq 0}}\Bigr\}
(-1)^i q^{\frac{1}{2}i(i+m)+p'j(pj+i)-\frac{1}{2}(\ell+1)(2pj+i)}.
\end{align*}

After these preliminaries let us now return to the conjugate Bailey pair 
of equation \eqref{gdinf} and specialize $a=q^{\eta}$, 
with $\eta$ a nonnegative integer.
Let us further remark the following 
identities ($\ell=0,1$ and $m+\ell\equiv L+\ell\equiv0\pmod{2}$):
$$C_{m,\ell}^{(1,3)}(q)=\frac{q^{\frac{1}{4}(m^2-\ell^2)}}{(q)_{\infty}}
\qquad \text{and} \qquad
X_{0,\ell+1}^{(1,3)}(L,1)=q^{\frac{1}{4}(L^2-\ell)}.$$
The first result is~\cite[Sec. 4.6, Ex. 3]{KP84} whereas the second is 
A(5) for $\ell=0$ and A(8) for $\ell=1$.
We thus infer that the conjugate Bailey pair \eqref{gdinf} can be recast
as $\gamma_L=(q)_{\eta}C_{2L+\eta,\ell}^{(1,3)}(q)$ and 
$\delta_L=X_{0,\ell+1}^{(1,3)}(2L+\eta,1)$.
It now requires little imagination to conjecture the following more general
result~\cite{SW99a}.
\begin{theorem}\label{thmCBP2}
Fix integers $1\leq p<p'$, and
let $\eta$ and $\ell$ be nonnegative integers such that $0\leq\ell\leq p'-2$ 
and $\ell+\eta$ is even. Then 
\begin{equation}\label{CBP2}
\gamma_L=(q)_{\eta}C_{2L+\eta,\ell}^{(p,p')}(q)  \qquad \text{and} \qquad
\delta_L=X^{(p,p')}_{0,\ell+1}(2L+\eta,1)
\end{equation}
yields a conjugate Bailey pair relative to $a=q^{\eta}$.
\end{theorem}
The proof of this theorem relies on yet another class of conjugate Bailey
pairs given by~\cite[Thm. 4.1]{SW99a}
\begin{equation}\label{K}
\begin{split}
\gamma_L&=\frac{1}{(q)_{\infty}^2(aq)_{\infty}}
\sum_{i=1}^{\infty}(-1)^i q^{\frac{1}{2}i(i+2L+\eta)}
\Bigl\{q^{\frac{1}{2}i(2j+\eta+1)}-q^{-\frac{1}{2}i(2j+\eta+1)}\Bigr\}\\[1mm]
\delta_L&=\qbins{2L+\eta}{L-j}-\qbins{2L+\eta}{L-j-1},
\end{split}
\end{equation}
with $a=q^{\eta}$, $\eta$ an nonnegative integer and $j$ an integer.
Here we remark that,
incidentally, $\delta_L=K_{(2^{L-j}1^{2j+\eta}),(1^{2L+\eta})}(q)=
\tilde{K}_{(L+j+\eta,L-j),(1^{2L+\eta})}(q)$, where $K_{\lambda,\mu}(q)$
and $\tilde{K}_{\lambda,\mu}(q)$ are the Kostka and cocharge
Kostka polynomial, respectively. The Bailey pair \eqref{K}
can easily be derived from the summation formula~\cite{SW99a}
\begin{equation*}
\sum_{r=0}^{\infty}\frac{q^r(ab)_{2r}}{(q)_r(ab)_r(aq)_r(bq)_r}
=\frac{1}{(q)_{\infty}(aq)_{\infty}(bq)_{\infty}}
\sum_{i=1}^{\infty}(-1)^{i-1} q^{\binom{i}{2}}\frac{a^i-b^i}{a-b}.
\end{equation*}

It is again possible to give representations of the polynomials
$X_{0,\ell+1}^{(p,p')}(L,1)$ that are manifestly positive~\cite{SW99a}.
Treating only the simplest cases we have the following
counterpart of Theorem~\ref{thmF}.
\begin{theorem}\label{thmF2}
Let $1<p<p'<2p$ with $\gcd(p,p')=1$. 
Then 
\begin{equation*}
X_{0,1}^{(p,p')}(2L,1)=q^L
\sum_{m\in 2\Z^d} q^{\frac{1}{4}mBm+\frac{1}{2}\sum_{i=1}^n m_{t_i}}
\prod_{j=1}^d\qbins{L\delta_{j,1}-\sum_{i=1}^n \delta_{j,t_i}
+\frac{1}{2}(\I m)_j}{m_j}
\end{equation*}
\end{theorem}
The corresponding identities for $2p<p'$ follow from
$$X_{r,s}^{(p'-p,p')}(2L,1;q) = q^{L(L+1)} X_{r,s}^{(p,p')}(2L,1;1/q).$$

Combining the Bailey pair of equation \eqref{alphabeta}
with the conjugate Bailey pair of \eqref{CBP2} and specializing
some of the parameters, we find
\begin{multline*}
\sum_{j\in\Z}\Bigl\{
q^{j(jp_1p_1'+rp_1'-sp_1)}C_{2jp_1',0}^{(p_2,p_2')}(q)
-q^{(jp_1+r)(jp_1'+s)}C_{2jp_1'+2s,0}^{(p_2,p_2')}(q)\Bigr\} \\ =
\sum_{L=0}^{\infty}X^{(p_1,p_1')}_{r,s}(2L,s)
X^{(p_2,p_2')}_{0,1}(2L,1)/(q)_{2L},
\end{multline*}
with $1\leq p_i<p_i'$ ($i=1,2$) and $|p_1'r-p_1 s|=1$.
Here we have used the symmetry $C_{m,\ell}^{(p,p')}=C_{-m,\ell}^{(p,p')}$.
Inserting the representations for the
one-dimensional configuration sums provided by 
Theorems~\ref{thmF} and \ref{thmF2} this turns into a class of
`rather' nontrivial $q$-series identities. For $p_1=1$ or $p_2=1$
the left-hand side of the above equation can be identified as a 
branching function of the coset pair
$(\text{A}_1^{(1)}\oplus\text{A}_1^{(1)},\text{A}_1^{(1)})$ at levels
$N_1=p_1/p_1-2$, $N_2=p_2'/p_2-2$ and $N_1+N_2$, respectively.

\section{Further reading}
To conclude our overview of half a century of Bailey's lemma,
let us mention some further papers on (or related to) the Bailey lemma
that have not been mentioned in the main text.
In \cite{Paule87}, Paule gave a short operator-type proof 
of the special case \eqref{iter2} of the Bailey chain.
Riese~\cite{Riese97} developed the Mathematica package \textsf{Bailey} for
taking (automated) walks along the Bailey lattice.
He also shows how to apply his Mathematica package \textsf{qZeil} to 
generate Bailey pairs.
The Bailey transform \eqref{BP} and its inverse
\eqref{inverseBP} can be formulated naturally in terms of 
inversion of infinite-dimensional lower-triangular 
matrices~\cite{AAB87,Bressoud88}, making it
a special case of the generalized $q$-Lagrange procedure of Gessel
and Stanton~\cite{GS83}.
New types of Baily lattice transformations which 
do not only change the value of $a$ but also that of the base $q$,
were very recently found  and applied by Bressoud, Ismail 
and Stanton~\cite{BIS99}.
Bressoud~\cite{Bressoud81} and Singh~\cite{Singh94} have also applied 
conjugate Bailey pairs other than \eqref{gd} and \eqref{gdinf} of Bailey.
For a special choice of parameters their conjugate pair can be shown
to coincide with the $(p,p')=(2,3)$ case of Theorem~\ref{thmCBP2}.
Andrews~\cite{Andrews86a}, Andrews and Hickerson~\cite{AH91} and 
Choi~\cite{Choi99} applied 
the Bailey chain to prove identities for Ramanujan's mock theta functions,
and Andrews~\cite{Andrews92} also used it to prove several
of Ramanujan's identities for Lambert series.
The Bailey lemma and its connection to $N=2$ supersymmetric conformal
field theory was investigated by Berkovich, McCoy and Schilling~\cite{BMS96}.
For those left with the impression that the Bailey lemma is
``merely'' good for proving $q$-series identities we remark that Andrews
utilized the Bailey machinery in \cite{Andrews86a,Andrews86b} to give
a proof of Gauss' theorem that every integer can be written as the sum
of three triangular numbers and that Andrews, Dyson and
Hickerson used Bailey's lemma in the context of algebraic number theory
\cite{ADH88}.
Finally we mention that a special case of the Bailey chain
admits an extension due to Burge~\cite{Burge93}. This was
extensively applied and further developed by Foda, Lee and Welsh~\cite{FLW98}
and by Schilling and the author~\cite{SW99b}.

\subsection*{Acknowledgements}
This work was supported by a fellowship of the Royal
Ne\-therlands Academy of Arts and Sciences.

\subsection*{Note added in proof}
The many recent references to the Bailey lemma listed in the 
bibliography show that after 50 years Bailey's lemma still is
a source of inspiration. This makes it quite impossible to 
publish an account that can claim to be complete and up to date.
Indeed, after this paper was accepted for publication further
advances in connection with the lemma were reported in 
\cite{Fulman00a,Fulman00b,Warnaar00}.

\end{document}